\documentstyle[12pt]{article}
\newcommand{\too}{\longrightarrow}

\newcommand{\di}{\displaystyle}
\newcommand{\Om}{\Omega}

\newcommand{\al}{\alpha}
\newcommand{\be}{\beta}
\newcommand{\ga}{\gamma}
\newcommand{\Ga}{\Gamma}

\newcommand{\De}{\Delta}

\def \reel{ {\rm I}\!{\rm R} }

 \def \rat{ {\rm Q}\kern-.65em {}^{{}_/ }}

\newtheorem{th}{Theorem}[section]
\newtheorem{pr}{Proposition}[section]

\title{Poisson structures compatible with the canonical metric of $\reel^3$} \author{M. Boucetta}
\date{} \parindent=0cm \begin{document}
\maketitle

{\bf Abstract.} In this Note, we will characterize the Poisson
structures compatible with the canonical metric of $\reel^3$. We
will also give some relevant examples of such structures. The
notion of compatibility between a Poisson structure and a
Riemannian metric used in this Note was introduced and studied by
the author in [1], [2], [3].

\section{Introduction and main results}
Many fundamental definitions and results about Poisson manifolds
can be found in Vaisman's monograph [5].

 As a continuation of  the study by
the author of  Poisson structures compatible with Riemannian
metrics in [1], [2] and [3], it is of interest to find some
relevant examples of such structures in the low dimensions. So we
were interested in finding all the Poisson structures compatible
with the canonical metric in  $\reel^3$. The results of this
search is the theme of this Note.

Let us recall some facts about the notion of compatibility between
a Poisson structure and a Riemannian metric in order to motivate
our investigation and to show the interest
  of this Note.

Let $P$ be a Poisson manifold with Poisson tensor $\pi$. A
Riemannian metric
  on  $T^*P$ is a smooth symmetric
contravariant 2-form $<,>$ on $P$ such that, at each point $x\in
P$, $<,>_x$ is a scalar product on $T_x^*P$.  For each Riemannian
metric $<,>$ on $T^*P$, we consider the contravariant connection
$D$ introduced in [1] by
\begin{eqnarray*} 2<D_\al\be,\ga>&=&\pi(\al).<\be
,\ga>+\pi(\be).<\al ,\ga>-\pi(\ga).< \al,\be>\\
&&+<[\al,\be]_{\pi},\ga >+<[\ga,\al]_{\pi}, \be>+<[\ga,\be]_{\pi},
\al>,\qquad(1)
\end{eqnarray*}where $\al,\be,\ga\in\Om^1(P)$ and the Lie bracket $[\;,\;]_\pi$
is given by $$
[\al,\be]_{\pi}=L_{\pi(\al)}\be-L_{\pi(\be)}\al-d(\pi(\al,\be));$$
here, $\pi:T^*P\too TP$ denotes the bundle map given by
$$\be[\pi(\al)]=\pi(\al,\be).$$  The connection $D$ is the
contravariant
 analogue of the usual  Levi-Civita
connection.
   The connection $D$ has vanishing torsion,
i.e.
 $$D_\al\be-D_\be\al=[\al,\be]_{\pi}.$$ Moreover, it is compatible with the
 Riemannian metric
 $<,>$,
 i.e.
 $$\pi(\al).<\be,\ga>=<D_\al\be,\ga>+<\be,D_\al\ga>.$$
 The notion of contravariant connection has been
introduced by Vaisman (see [5] p.55) as contravariant derivative.
Recently,
  a geometric approach of this notion was given by Fernandes in
  [4].

If we put, for any $f\in C^\infty(P)$,
 $$\phi_{<,>}(f)=
\sum_{i=1}^n<D_{\al_i}df,\al_i>\eqno(2)$$where
$(\al_1,\ldots,\al_n)$ is a local othonormal basis of 1-forms, we
get a derivation on $C^\infty(P)$ and hence  a vector field called
the modular vector field of $(P,\pi)$ with respect to the metric
$<,>$.

 The couple $(\pi,<,>)$ is compatible if, for any $\al,\be,\ga\in\Om^1(P)$,
$$D\pi(\al,\be,\ga):=
\pi(\al).\pi(\be,\ga)-\pi(D_\al\be,\ga)-\pi(\be,D_\al\ga)=0.\eqno(3)$$
In this case, the triple $(P,\pi,<,>)$ is called a Riemann-Poisson
manifold.

  Riemann-Poisson manifolds was first introduced by the
author in [1]. Let us summarize some important  results of
Riemann-Poisson manifolds proved by the author in [2] and [3].

For a Riemann-Poisson manifold $(P,\pi,<,>)$ the following results
are true:

1. the symplectic leaves are K\"ahlerian;

2. the symplectic foliation (when it is a regular foliation) is a
Riemannian foliation;

3. $(P,\pi)$ is unimodular (see [6] for the details on the notion
of unimodular Poisson manifolds) and moreover the modular vector
field $\phi_{<,>}$ given by $(2)$ vanishes.

With those properties in mind, we can give the main results of
this Note.
\begin{th} A Poisson tensor
$\pi=\pi_{12}\frac{\partial}{\partial x
}\wedge\frac{\partial}{\partial
y}+\pi_{13}\frac{\partial}{\partial
x}\wedge\frac{\partial}{\partial
z}+\pi_{23}\frac{\partial}{\partial
y}\wedge\frac{\partial}{\partial z}$  is compatible with the
canonical metric $<,>$ of $\reel^3$ iff there exists a function
$f\in C^\infty(\reel^3)$ such that
$$\pi_{12}=\frac{\partial f}{\partial z},\qquad\pi_{13}=-
\frac{\partial f}{\partial y},\qquad \pi_{23}=\frac{\partial
f}{\partial x},$$and
$$d(<df,df>)-\De(f)df=0\eqno(E)$$where $\di\De=\frac{\partial^2}{\partial x^2}+
\frac{\partial^2}{\partial y^2}+\frac{\partial^2}{\partial z^2}$
is the usual Laplacian on $\reel^3$. Moreover, the function $f$ is
a Casimir function of $\pi$.\end{th} The following proposition and
the theorem above give all the linear Poisson structures on
$\reel^3$ compatible with the canonical metric.
\begin{pr} The polynomial functions of degree 2 solutions of $(E)$
are
$$f(x,y,z)=(a+c)x^2+(a+b)y^2+(b+c)z^2-2\sqrt{bc}xy+2\sqrt{ab}xz+
2\sqrt{ac}yz,$$where $a,b,c\in\reel$ and
$ab,ac,bc\in\reel_+$.\end{pr}

Let $\pi_{so(3)}=z\frac{\partial}{\partial x
}\wedge\frac{\partial}{\partial y}-y\frac{\partial}{\partial
x}\wedge\frac{\partial}{\partial z}+x\frac{\partial}{\partial
y}\wedge\frac{\partial}{\partial z}$ be the linear Poisson
structure on $\reel^3$ corresponding to the Lie algebra $so(3)$.
In $[2]$, we have shown that there isn't any Riemannian metric on
$\reel^3$ compatible with $\pi_{so(3)}$. However, we have the
following  proposition.
\begin{pr} The function $f(x,y,z)=(x^2+y^2+z^2)^{\frac32}$ is a
solution of $(E)$ and then $(x^2+y^2+z^2)^{\frac12}\pi_{so(3)}$ is
compatible with the canonical metric of $\reel^3$.\end{pr}

{\bf Remarks.} 1. The fact that there isn't any metric compatible
with $\pi_{so(3)}$ and, however, the Poisson structure
$(x^2+y^2+z^2)^{\frac12}\pi_{so(3)}$ is compatible with the
canonical metric seems curious. But it can be explained easily. In
fact, let $(P,\pi,<,>)$ be a Poisson manifold with a contravariant
Riemannian metric. If we change the Poisson structure by $f\pi$
where $f\in C^\infty(P)$, the contravariant Levi-Civita connection
given by (1) become more complicated and is given by
$$D^{f\pi}_{\al}\be=fD^{\pi}_{\al}\be+\frac12\pi(\al,\be)df-\frac12<df,\be>J\al
-\frac12<df,\al>J\be$$where $J$ is the field of homomorphisms
given by $\pi(\al,\be)=<J\al,\be>$.

2. The Poisson structures $\pi_{so(3)}$ and
$(x^2+y^2+z^2)^{\frac12}\pi_{so(3)}$ have the same symplectic
foliation and, in restriction to a symplectic leaf, the two
symplectic structures differ by a constant.

3. It is possible that the compatibility of
$(x^2+y^2+z^2)^{\frac12}\pi_{so(3)}$ with the canonical metric has
some physical signification.
\section{Proof of Theorem 1.1}
Let  $\pi$ be a bivectors  field on $\reel^3$ given by
$$\pi=\pi_{12}\frac{\partial }{\partial x}\wedge\frac{\partial }{\partial y}
+\pi_{13}\frac{\partial }{\partial x}\wedge\frac{\partial
}{\partial z}+\pi_{23}\frac{\partial }{\partial
y}\wedge\frac{\partial }{\partial z}.$$ We consider the canonical
metric $<,>$ on $\reel^3$ as contravariant metric given by
$$<dx,dx>=<dy,dy>=<dz,dz>=1, <dx,dy>=<dx,dz>=<dy,dz>=0.$$
We denote by $D$ the Levi-Civita contravariant connection
associated with $(\pi,<,>)$.

 Firstly, remark that the compatibility
between $\pi$ and $<,>$ implies the vanishing of the modular
vector field given by
$$\phi_{<,>}(f)=<D_{dx}df,dx>+<D_{dy}df,dy>+<D_{dz}df,dz>.$$
A straightforward  calculation shows that the vanishing of
$\phi_{<,>}$ is equivalent to
$$
\frac{\partial\pi_{12}}{\partial
y}+\frac{\partial\pi_{13}}{\partial z}=0,\qquad
\frac{\partial\pi_{12}}{\partial
x}-\frac{\partial\pi_{23}}{\partial z}=0,\qquad
\frac{\partial\pi_{13}}{\partial
x}+\frac{\partial\pi_{23}}{\partial y}=0.\eqno(4)$$ Now, it is
easy to see that $(4)$ is equivalent to the fact that
$\pi_{23}dx-\pi_{13}dy+\pi_{12}dz$ is a closed 1-form and hence is
exact. So, there exists a function $f\in C^\infty(\reel^3)$ such
that
$$\pi_{12}=\frac{\partial f}{\partial z},\qquad\pi_{13}=-
\frac{\partial f}{\partial y},\qquad \pi_{23}=\frac{\partial
f}{\partial x}.\eqno(5)$$ Now, let us compute the contravariant
connection $D$. We will use the Christoffel symbols $\Ga_{ij}^k$.
For example, $D_{dx}dx=\Ga_{11}^1dx+\Ga_{11}^2dy+\Ga_{11}^3dz$.
From (1), we  get:
\begin{eqnarray*}
&&\Ga_{11}^1=0,\quad\Ga_{11}^2=-\frac{\partial^2f}{\partial x
\partial z},\quad
\Ga_{11}^3=\frac{\partial^2f}{\partial x\partial y},\\
&&\Ga_{12}^1=\frac{\partial^2f}{\partial x\partial
z},\quad\Ga_{12}^2=0,\quad
\Ga_{12}^3=\frac12\left(-\frac{\partial^2f}{\partial
x^2}+\frac{\partial^2f}{\partial y^2}
+\frac{\partial^2f}{\partial z^2}\right),\\
&&\Ga_{21}^1=0,\quad\Ga_{21}^2=-\frac{\partial^2f}{\partial
y\partial z},\quad
\Ga_{21}^3=\frac12\left(-\frac{\partial^2f}{\partial
x^2}+\frac{\partial^2f}{\partial y^2}
-\frac{\partial^2f}{\partial z^2}\right),\\
 &&\Ga_{13}^1=-\frac{\partial^2f}{\partial x\partial
y},\quad\Ga_{13}^2=\frac12\left(\frac{\partial^2f}{\partial
x^2}-\frac{\partial^2f}{\partial y^2} -\frac{\partial^2f}{\partial
z^2}\right),\quad
\Ga_{13}^3=0,\\
 &&\Ga_{31}^1=0,\quad\Ga_{31}^2=\frac12\left(\frac{\partial^2f}{\partial
x^2}+\frac{\partial^2f}{\partial y^2} -\frac{\partial^2f}{\partial
z^2}\right),\quad\Ga_{31}^3=\frac{\partial^2f}{\partial y\partial
z},\\
&&\Ga_{22}^1=\frac{\partial^2f}{\partial y\partial
z},\quad\Ga_{22}^2=0,\quad
\Ga_{22}^3=-\frac{\partial^2f}{\partial x\partial y},\\
&&\Ga_{23}^1=\frac12\left(\frac{\partial^2f}{\partial
x^2}-\frac{\partial^2f}{\partial y^2} +\frac{\partial^2f}{\partial
z^2}\right),\quad\Ga_{23}^2=\frac{\partial^2f}{\partial x\partial
y},\quad
\Ga_{23}^3=0,\\
&&\Ga_{32}^1=\frac12\left(-\frac{\partial^2f}{\partial
x^2}-\frac{\partial^2f}{\partial y^2} +\frac{\partial^2f}{\partial
z^2}\right),\quad\Ga_{32}^2=0,\quad\Ga_{32}^3=-\frac{\partial^2f}{\partial
x\partial z},\\
&&\Ga_{33}^1=-\frac{\partial^2f}{\partial y\partial
z},\quad\Ga_{33}^2=\frac{\partial^2f}{\partial x\partial z}
,\quad\Ga_{33}^3=0.\\
\end{eqnarray*}
Now, we will compute $D_{dx}\pi$, $D_{dy}\pi$ and $D_{dz}\pi$. We
have
\begin{eqnarray*}
D_{dx}\pi&=&\pi(dx)(\frac{\partial f}{\partial
z})\frac{\partial}{\partial x }\wedge\frac{\partial}{\partial
y}-\pi(dx)(\frac{\partial f}{\partial y})\frac{\partial}{\partial
x }\wedge\frac{\partial}{\partial z}+\pi(dx)(\frac{\partial
f}{\partial x})\frac{\partial}{\partial y
}\wedge\frac{\partial}{\partial z}\\
&+&\frac{\partial f}{\partial
z}\left((D_{dx}\frac{\partial}{\partial x
})\wedge\frac{\partial}{\partial y}+\frac{\partial}{\partial x
}\wedge(D_{dx}\frac{\partial}{\partial y})\right)\\
&-&\frac{\partial f}{\partial
y}\left((D_{dx}\frac{\partial}{\partial x
})\wedge\frac{\partial}{\partial z}+\frac{\partial}{\partial x
}\wedge(D_{dx}\frac{\partial}{\partial z})\right)\\
&+&\frac{\partial f}{\partial
x}\left((D_{dx}\frac{\partial}{\partial y
})\wedge\frac{\partial}{\partial z}+\frac{\partial}{\partial y
}\wedge(D_{dx}\frac{\partial}{\partial z})\right).
\end{eqnarray*}
On other hand, we have
\begin{eqnarray*}
\pi(dx)&=&\frac{\partial f}{\partial z}\frac{\partial}{\partial
y}-\frac{\partial f}{\partial y}\frac{\partial}{\partial z},\\
D_{dx}\frac{\partial}{\partial x }&=&-\frac{\partial^2f}{\partial
x\partial z}\frac{\partial}{\partial y}
+\frac{\partial^2f}{\partial x\partial y}\frac{\partial}{\partial
z},\\
D_{dx}\frac{\partial}{\partial y }&=&\frac{\partial^2f}{\partial
x\partial z}\frac{\partial}{\partial x}
+\frac12\left(-\frac{\partial^2f}{\partial
x^2}+\frac{\partial^2f}{\partial y^2} +\frac{\partial^2f}{\partial
z^2}\right)\frac{\partial}{\partial z},\\
D_{dx}\frac{\partial}{\partial z }&=&-\frac{\partial^2f}{\partial
x\partial y}\frac{\partial}{\partial x}
+\frac12\left(\frac{\partial^2f}{\partial
x^2}-\frac{\partial^2f}{\partial y^2} -\frac{\partial^2f}{\partial
z^2}\right)\frac{\partial}{\partial y}.\end{eqnarray*}
Substituting those expressions into the expression of $D_{dx}\pi$,
we get
\begin{eqnarray*}
D_{dx}\pi&=&\left(\frac{\partial f}{\partial
z}\frac{\partial^2f}{\partial y\partial z}+\frac{\partial
f}{\partial x}\frac{\partial^2f}{\partial x\partial
y}+\frac12\frac{\partial f}{\partial
y}\left(-\frac{\partial^2f}{\partial
x^2}+\frac{\partial^2f}{\partial y^2} -\frac{\partial^2f}{\partial
z^2}\right)\right)\frac{\partial}{\partial x
}\wedge\frac{\partial}{\partial y}\\
&+&\left(\frac{\partial f}{\partial y}\frac{\partial^2f}{\partial
y\partial z}+\frac{\partial f}{\partial
x}\frac{\partial^2f}{\partial x\partial z}+\frac12\frac{\partial
f}{\partial z}\left(-\frac{\partial^2f}{\partial
x^2}-\frac{\partial^2f}{\partial y^2} +\frac{\partial^2f}{\partial
z^2}\right)\right)\frac{\partial}{\partial x
}\wedge\frac{\partial}{\partial z}.\end{eqnarray*} In the same
manner we can get
\begin{eqnarray*}
D_{dy}\pi&=&\left(-\frac{\partial f}{\partial
z}\frac{\partial^2f}{\partial x\partial z}-\frac{\partial
f}{\partial y}\frac{\partial^2f}{\partial x\partial
y}+\frac12\frac{\partial f}{\partial
x}\left(-\frac{\partial^2f}{\partial
x^2}+\frac{\partial^2f}{\partial y^2} +\frac{\partial^2f}{\partial
z^2}\right)\right)\frac{\partial}{\partial x
}\wedge\frac{\partial}{\partial y}\\
&+&\left(\frac{\partial f}{\partial y}\frac{\partial^2f}{\partial
y\partial z}+\frac{\partial f}{\partial
x}\frac{\partial^2f}{\partial x\partial z}+\frac12\frac{\partial
f}{\partial z}\left(-\frac{\partial^2f}{\partial
x^2}-\frac{\partial^2f}{\partial y^2} +\frac{\partial^2f}{\partial
z^2}\right)\right)\frac{\partial}{\partial y
}\wedge\frac{\partial}{\partial z}.\\
D_{dz}\pi&=&\left(-\frac{\partial f}{\partial
z}\frac{\partial^2f}{\partial x\partial z}-\frac{\partial
f}{\partial y}\frac{\partial^2f}{\partial x\partial
y}+\frac12\frac{\partial f}{\partial
x}\left(-\frac{\partial^2f}{\partial
x^2}+\frac{\partial^2f}{\partial y^2} +\frac{\partial^2f}{\partial
z^2}\right)\right)\frac{\partial}{\partial x
}\wedge\frac{\partial}{\partial z}\\
&+&\left(-\frac{\partial f}{\partial z}\frac{\partial^2f}{\partial
y\partial z}-\frac{\partial f}{\partial
x}\frac{\partial^2f}{\partial x\partial y}+\frac12\frac{\partial
f}{\partial y}\left(+\frac{\partial^2f}{\partial
x^2}-\frac{\partial^2f}{\partial y^2} +\frac{\partial^2f}{\partial
z^2}\right)\right)\frac{\partial}{\partial y
}\wedge\frac{\partial}{\partial z}.
\end{eqnarray*}
Now, it is easy to show that $D\pi=0$ iff $f$ satisfies $(E)$. It
is also easy to show that $f$ is a Casimir function. Remark that
 $D\pi=0$ implies that the bracket of Schouten $[\pi,\pi]$
vanishes which finish the proof of Theorem 1.1.$\Box$

{\bf References}

\bigskip

[1] M. Boucetta,  Compatibilit\'e des structures
pseudo-riemanniennes et des structures de Poisson, C. R. Acad.
Sci. Paris, {\bf t. 333}, S\'erie I, (2001) 763-768.

[2] M. Boucetta, Poisson manifolds with compatible pseudo-metric
and pseudo-Riemannian Lie algebras, Preprint math.DG/0206102. To
appear in  Differential Geometry and its Applications.

[3] M. Boucetta, Riemann-Poisson manifolds and K\"ahler-Riemann
foliations, C. R. Acad. Sci. Paris, Ser. I 336 (2003) 423-428.

[4]  R. L. Fernandes, Connections in Poisson Geometry 1: holonomy
and invariants, J. Diff. Geom. {\bf 54} , (2000) 303-366.

[5] I. Vaisman, Lectures on the Geometry of Poisson Manifolds,
Progress in Mathematics, vol. {\bf118}, Birkh\"auser, Berlin,
1994.

[6] A. Weinstein, The Modular Automorphism Group of a Poisson
Manifold, J. Geom. Phys. {\bf23}, (1997) 379-394.

\bigskip
{\it Mohamed Boucetta

Facult\'e des Sciences et Techniques Gueliz

BP 549 Marrakech Morocco

E-mail: boucetta@fstg-marrakech.ac.ma}

\end{document}